\documentclass{amsart}

\usepackage[all]{xy}
\usepackage{amssymb}
\usepackage{amsmath}
\usepackage{graphicx}
\usepackage[english]{babel}
\usepackage{color}
\usepackage{enumerate}

\newcommand{\comment}[1]{}

\theoremstyle{definition}
\newtheorem{definition}{Definition}
\newtheorem{remark}[definition]{Remark}
\newtheorem{example}[definition]{Example}

\theoremstyle{plain}
\newtheorem{lemma}[definition]{Lemma}
\newtheorem{proposition}[definition]{Proposition}
\newtheorem{corollary}[definition]{Corollary}
\newtheorem{theorem}[definition]{Theorem}

\def\NN{{\mathbb{N}}_{\geq 2}}
\def\NNN{{\mathbb{N}}}
\def\KK{{\mathbb{K}}}

\def\HH{{\mathcal{H}}}

\def\TT{{\mathcal{T}}}
\def\PP{{\mathcal{P}}}

\def\MM{{\mathcal{M}}}
\def\RRR{{\mathcal{R}}}
\def\III{{\mathcal{I}}}
\def\QQQ{{\mathcal{Q}}}

\def\Prim{\mathrm{Prim}}
\def\Id{\mathrm{Id }}

\def\PT{\mathrm{PRT}}
\def\root{\mathit{root}}

\def\Sh{\mathop{\rm Sh}}

\def\m{\underline m}
\def\n{\underline n}

\def\Mag{\mathcal Mag}
\def\Nil{\mathcal Nil}
\def\Vect{\mathcal Vect}

\newcommand{\mbialg}[2]{$(#1; #2)$-magmatic bialgebra}
\newcommand{\malg}[1]{$#1$-magmatic algebra}
\newcommand{\mop}[1]{$#1$-magmatic operad}
\newcommand{\mcoalg}[1]{$#1$-magmatic coalgebra}
\newcommand{\mrootalg}[2]{$#2_{#1}$-rooted magmatic algebra}
\newcommand{\mrootop}[2]{$#2_{#1}$-rooted magmatic operad}
\newcommand\magroot[2]{\mathcal{M}ag^{#2}_{\root\setminus #1}}

\newcommand{\F}[1]{F_#1\ }

\def\Prim{\operatorname{Prim\ }}

\def\Id{\mathrm{Id}}
\def\Lie{\mathcal{L}ie}

\def\arbreA {\vcenter{\xymatrix@R=3pt@C=3pt{
&& \\
&*{}\ar@{-}[ur] \ar@{-}[ul] \ar@{-}[d]     &\\
&&}}}

\def\arbreAt {\vcenter{\xymatrix@R=3pt@C=3pt{
&& \\
&& \\
&*{}\ar@{-}[uur] \ar@{-}[uu]\ar@{-}[uul] \ar@{-}[d]     &\\
&&}}}

\def\arbreAttt{\vcenter{\xymatrix@R=3pt@C=3pt{
&&&& \\
&&&&\\
&&*{}\ar@{-}[uurr] \ar@{-}[uur]\ar@{-}[uul] \ar@{-}[uull]\ar@{-}[dd]     &&\\
&&&&\\
&&&&}}}

\def\mcorol{\vcenter{\xymatrix@R=3pt@C=3pt{
&&m&& \\
&&\cdots&&\\
&&*{}\ar@{-}[uurr] \ar@{-}[uur]\ar@{-}[uul] \ar@{-}[uull]\ar@{-}[d]     &&\\
&&&&}}}

\def\arbreBA{\vcenter{\xymatrix@R=2pt@C=2pt{
&&&&\\
&&&*{}\ar@{-}[ul] & \\
&&*{}\ar@{-}[uurr] \ar@{-}[uull] \ar@{-}[d]     &&\\
&&&&}}}

\def\arbreBAt{\vcenter{\xymatrix@R=2pt@C=2pt{
&&&&\\
&&&*{}\ar@{-}[ul] & \\
&&*{}\ar@{-}[uurr] \ar@{-}[uu]\ar@{-}[uull] \ar@{-}[dd]     &&\\
&&&&\\
&&&&}}}

\def\arbreBAtt{\vcenter{\xymatrix@R=2pt@C=2pt{
&&&&\\
&&&*{}\ar@{-}[ul]\ar@{-}[u] & \\
&&*{}\ar@{-}[uurr] \ar@{-}[uu]\ar@{-}[uull] \ar@{-}[d]     &&\\
&&&&}}}

\def\arbreBAttt{\vcenter{\xymatrix@R=2pt@C=2pt{
&&&&\\
&&&*{}\ar@{-}[ul]\ar@{-}[u] & \\
&&*{}\ar@{-}[uurr] \ar@{-}[uull] \ar@{-}[dd]     &&\\
&&&&\\
&&&&}}}

\def\arbreAB{\vcenter{\xymatrix@R=2pt@C=2pt{
&&&&\\
&*{}\ar@{-}[ur] &&& \\
&&*{}\ar@{-}[uurr] \ar@{-}[uull] \ar@{-}[d]     &&\\
&&&&}}}

\def\arbreABttt{\vcenter{\xymatrix@R=2pt@C=2pt{
&&&&\\
&*{}\ar@{-}[ur]\ar@{-}[u] &&& \\
&&*{}\ar@{-}[uurr] \ar@{-}[uull] \ar@{-}[dd]     &&\\
&&&&\\
&&&&}}}

\def\arbreABt{\vcenter{\xymatrix@R=2pt@C=2pt{
&&&&\\
&*{}\ar@{-}[ur] &&& \\
&&*{}\ar@{-}[uurr] \ar@{-}[uu]\ar@{-}[uull] \ar@{-}[dd]     &&\\
&&&&\\
&&&&}}}

\def\arbreABtt{\vcenter{\xymatrix@R=2pt@C=2pt{
&&&&\\
&*{}\ar@{-}[ur] \ar@{-}[u]&&& \\
&&*{}\ar@{-}[uurr] \ar@{-}[uu]\ar@{-}[uull] \ar@{-}[dd]     &&\\
&&&&\\
&&&&}}}

\def\arbreAAt{\vcenter{\xymatrix@R=2pt@C=2pt{
&&&&\\
&&*{}\ar@{-}[ur]\ar@{-}[ul] && \\
&&*{}\ar@{-}[uurr] \ar@{-}[u]\ar@{-}[uull] \ar@{-}[dd]     &&\\
&&&&\\
&&&&}}}

\def\arbreAAtt{\vcenter{\xymatrix@R=2pt@C=2pt{
&&&&\\
&&*{}\ar@{-}[ur] \ar@{-}[u]&& \\
&&*{}\ar@{-}[uurr] \ar@{-}[uu]\ar@{-}[uull] \ar@{-}[dd]     &&\\
&&&&\\
&&&&}}}

\def\arbreBB{\vcenter{\xymatrix@R=2pt@C=2pt{
&&&&\\
&&&& \\
&&*{}\ar@{-}[uurr] \ar@{-}[uull] \ar@{-}[d] \ar@{-}[uu]     &&\\
&&&&}}}

\def\arbreABC{\vcenter{\xymatrix@R=1pt@C=1pt{
&&&&&&\\
&*{}\ar@{-}[ur] &&&&& \\
&&*{}\ar@{-}[uurr] &&&&\\
&&&*{}\ar@{-}[uuurrr] \ar@{-}[uuulll] \ar@{-}[d] &&&\\
&&&&&&}}}

\def\arbreBAC{\vcenter{\xymatrix@R=1pt@C=1pt{
&&&&&&\\
&&&*{}\ar@{-}[ul] &&& \\
&&*{}\ar@{-}[uurr] &&&&\\
&&&*{}\ar@{-}[uuurrr] \ar@{-}[uuulll] \ar@{-}[d] &&&\\
&&&&&&}}}

\def\arbreACA{\vcenter{\xymatrix@R=1pt@C=1pt{
&&&&&&\\
&*{}\ar@{-}[ur] &&&&*{}\ar@{-}[ul] & \\
&&&&&&\\
&&&*{}\ar@{-}[uuurrr] \ar@{-}[uuulll] \ar@{-}[d] &&&\\
&&&&&&}}}

\def\arbreCAB{\vcenter{\xymatrix@R=1pt@C=1pt{
&&&&&&\\
&&&*{}\ar@{-}[ur] &&& \\
&&&&*{}\ar@{-}[uull] &&\\
&&&*{}\ar@{-}[uuurrr] \ar@{-}[uuulll] \ar@{-}[d] &&&\\
&&&&&&}}}

\def\arbreCBA{\vcenter{\xymatrix@R=1pt@C=1pt{
&&&&&&\\
&&&&&*{}\ar@{-}[ul] & \\
&&&&*{}\ar@{-}[uull] &&\\
&&&*{}\ar@{-}[uuurrr] \ar@{-}[uuulll] \ar@{-}[d] &&&\\
&&&&&&}}}

\def\arbreCBAt{\vcenter{\xymatrix@R=1pt@C=1pt{
&&&&&&\\
&&&&&*{}\ar@{-}[ul]\ar@{-}[u] & \\
&&&&*{}\ar@{-}[uull]\ar@{-}[uu] &&\\
&&&*{}\ar@{-}[uuurrr] \ar@{-}[uuulll] \ar@{-}[d] &&&\\
&&&&&&}}}

\def\stage{\vcenter{\xymatrix@R=1pt@C=1pt{
&&&&\\
&*{}\ar@{-}[ul]&\cdots &*{}\ar@{-}[ur]&n\\
*{}\ar@{--}[rrrr]&&&&\\
&&\cdots&&\\
*{}\ar@{--}[rrrr]&&&&\\
&&&& \\
&&*{}\ar@{-}[ur] \ar@{-}[ul] \ar@{-}[d]     &&\\
&&&&1 \ .
}}}

\def\concattt{\vcenter{\xymatrix@R=1pt@C=1pt{
t_1&&t_2\\
*{}\ar@{-}[u]&&*{}\ar@{-}[u]\\
&*{}\ar@{-}[ul] \ar@{-}[dd] \ar@{-}[ur]&\\
&&\\
&&}}}

\def\concat{\vcenter{\xymatrix@R=1pt@C=1pt{
t&&s\\
*{}\ar@{-}[u]&&*{}\ar@{-}[u]\\
&*{}\ar@{-}[ul] \ar@{-}[dd] \ar@{-}[ur]&\\
&&\\
&&}}}

\def\concatt{\vcenter{\xymatrix@R=1pt@C=1pt{
t&s&u\\
*{}\ar@{-}[u]&&*{}\ar@{-}[u]\\
&*{}\ar@{-}[ul] \ar@{-}[uu]\ar@{-}[dd] \ar@{-}[ur]&\\
&&\\
&&}}}

\def\deconcattwo{\vcenter{\xymatrix@R=1pt@C=1pt{
t_1&&t_2\\
&\otimes&\\
&&\\
*{}\ar@{-}[uuu]&&*{}\ar@{-}[uuu]\\
}}}

\def\deconcat{\vcenter{\xymatrix@R=1pt@C=1pt{
t_1&&&&t_n\\&\otimes&\cdots&\otimes&\\
&&&&\\
*{}\ar@{-}[uuu]&&&&*{}\ar@{-}[uuu]\\
}}}

\def\deconcatk{\vcenter{\xymatrix@R=1pt@C=1pt{
t_1&&&&t_k\\&\otimes&\cdots&\otimes&\\ &&&&\\
*{}\ar@{-}[uuu]&&&&*{}\ar@{-}[uuu]\\ }}}

\def\arbreetiquette{\xymatrix@R=1pt@C=1pt{
v_{1}\ar@{-}[dd] &v_{2}\ar@{-}[dd] &...&&v_{n}\ar@{-}[dd] \\
&&&&&\\
&&&&&\\
&&...&&\\
&&&&&\\
&&&&&\\
&&*{}\ar@{-}[uurr] \ar@{-}[uull] \ar@{-}[d] &&&\\
&&&&& \ .
}}

\def\colibre{\xymatrix@R=1pt@C=1pt{
&&&&&&&\\
C \ar@{>}[ddddrrrrrrr]_{\phi} \ar@{-->}[rrrrrr]^{\tilde \phi}&&&&&&&C(V)^\infty=\oplus_{n \geq 0} C_{n}^\infty\otimes V^{\otimes n} \ar@{>>}[dddd]\\
&&&&&&&\\
&&&&&&&\\
&&&&&&&\\
&&&&&&&V \ .
}}

\def\concatcolor{\xymatrix@R=1pt@C=1pt{
T&&W\\
*{}\ar@{-}[u]&&*{}\ar@{-}[u]\\
&x{}\ar@{-}[ul] \ar@{-}[dd] \ar@{-}[ur]& \\
&&\\
&&}}

\def\produit{\xymatrix@R=1pt@C=1pt{
A^{\otimes k}\ar@{->}[dddd]_{\mu_k}&&&&&&&&\ar@{->}[llllllll]_-{\Id\otimes\cdots\otimes \eta\otimes\cdots\otimes\Id}A^{\otimes i}\otimes \mathbb K \otimes A^{\otimes k-i-1}\ar@{=}[rr]&&\ar@{->}[ddddllllllllll]^-{\mu_{k-1}}A^{\otimes k-1}\\
&&&&&&&&\\
&&&&&&&&\\
&&&&&&&&\\
A &&&&&&&&\\
&&&&&&&&&&&\ }}

\def\coproduit{\xymatrix@R=1pt@C=1pt{
C^{\otimes k}\ar@{->}[rrrrrrrr]^-{\Id\otimes\cdots\otimes \epsilon\otimes\cdots\otimes\Id}&&&&&&&&C^{\otimes i}\otimes \mathbb K \otimes C^{\otimes k-i-1}\ar@{=}[rr]&&C^{\otimes k-1}\\
&&&&&&&&\\
&&&&&&&&\\
&&&&&&&\\
C\ar@{->}[uuuu]_{\Delta_k}\ar@{->}[uuuurrrrrrrrrr]_{\Delta_{k-1}}&&&&&&&\\
&&&&&&&&&&&
}}

\begin{document}
\author[E. Burgunder, R. Holtkamp]{Emily Burgunder, Ralf Holtkamp}
\address{EB:
Institut de Math\'ematiques et de mod\'elisation de Montpellier \\
UMR CNRS 5149\\
D\'epartement de math\'ematiques\\
Universit\'e Montpellier II\\
Place Eug\`ene Bataillon\\
34095 Montpellier CEDEX 5\\
France}
\email{burgunder@math.univ-montp2.fr}
\address{RH: Fakult\"at f\"ur Mathematik, Ruhr-Universit\"at,
          44780 Bochum, Germany}
\email{ralf.holtkamp@ruhr-uni-bochum.de}

\title{Partial magmatic bialgebras}

\date{\today}

\begin{abstract}
A partial magmatic bialgebra, \mbialg T S where $T\subset S$ are
subsets of $\NN$, is a vector space endowed with an $n$-ary
operation for each $n\in S$ and an $m$-ary co-operation for each
$m\in T$ satisfying some compatibility and unitary relations. We
prove an analogue of the Poincar\'e-Birkhoff-Witt theorem for
these partial magmatic bialgebras.
\end{abstract}

\maketitle

\section{Introduction}

Let $\HH$ be a cocommutative bialgebra over a field $\KK$ of
characteristic zero. Then the Poincar\'e-Birkhoff-Witt and the
Cartier-Milnor-Moore theorem \cite{Ca,M} can be rephrased as:

\begin{theorem}[PBW CMM]
 The following are equivalent:
\begin{enumerate}[i.]
\item $\HH$ is connected, \item $\HH$ is isomorphic to
$U(\Prim\HH)$, \item $\HH$ is cofree among the connected  cofree
coalgebras.
\end{enumerate}
The space of primitive elements $\Prim\HH$ of $\HH$ functorially
admits a $\Lie$-algebra structure, and $U$ is the universal
enveloping algebra functor.
\end{theorem}

\medskip

This structure theorem has been extended to connected
associative bialgebras by J.-L. Loday and M. Ronco in
\cite{LR1}, and has lead to many other structure theorems for
generalised bialgebras (see for example \cite{H,HLR,M,Q,R1,R}).

There exists now a general theory for these structure theorems for
non-unital bialgebras over operads due to J.-L. Loday, see
\cite{L}. These bialgebra types are given by  operads  ${\mathcal
A}$ handling the algebra structure,  operads ${\mathcal C}$
handling the coalgebra structure, and the compatibility relations
between the two structures. This theory gives conditions for the
existence of a primitive operad such that an analogue of the
structure theorem holds. Even when the conditions are verified,
one of the main difficulties remains in unravelling the algebraic
structure of the primitives.

Motivated by a systematic study of the theory, the first
paradigmatic example of a bialgebra type to consider is the case
where both operads ${\mathcal A}$ and ${\mathcal C}$ are free. We
call this bialgebra type the partial magmatic bialgebras, and
study them in this article, continuing some works of both authors
\cite{B1,B2,H,HLR}, which consider certain bialgebra types with at
least one free operad $\Mag$ -- or $\Mag^{\underline n}$, with
free generators up to arity $n$. Magmatic bialgebras defined over
${\mathcal A}={\mathcal C}=\Mag^{\underline n}$ were introduced
and classified - over the operad ${\mathcal Vect}$ of primitives -
by the first author in \cite{B2}. The natural generalisation of
the case of magmatic bialgebras over ${\mathcal
C}=\Mag^{\underline m}$, ${\mathcal A}=\Mag^{\underline n}$, for
$m\neq n$, was posed as an open problem in \cite{B2}.

More generally, we consider in the first part of the paper \mbialg
T Ss where $T\subset S\subset\NNN$. Such a \mbialg T S is a vector
space endowed with several not necessarily (co)associative
operations and co-operations, namely one $n$-ary operation for
each $n\in S$ and one $m$-ary co-operation for each $m\in T$. The
given bialgebra type verifies the conditions of \cite{L} if the
operations and cooperations are supposed to verify a compatibility
relation, namely the partial magmatic compatibility relation.  The
primitive operads $\magroot T S$ that we unravel are free operads
generated by infinitely many generators. An equivalent way of
defining the operad $\magroot T S$ is as $\KK\ \Id \oplus\III$,
where $\III$ is a right operadic ideal of the operad $\Mag^S$. The
structure theorem holds for \mbialg T Ss, with $\Lie$ replaced by
$\magroot T S$.

The second part of the article is to adapt this theorem to the
unital framework. In order to stay in the operadic framework, we have to restrict our study to partial
magmatic bialgebras where $S=\{2,\cdots,n\}$ and
$T=\{2,\cdots,m\}$.

The question is how the statements of the first part have to be
modified so that they are  still valid in the presence of a unit.
We give an answer to this question by finding the precise
compatibility relation that one has to take, and by the way answer
the open question of \cite{B2}.

Since the composition of operations may no longer be associative,
the unital framework for general partial magmatic bialgebras needs
a more general setting than the operadic setting.

\medskip

In Section 2 we recall the basics on planar trees, and we
introduce \malg Ss, \mcoalg Ts. From Section 3 on, we work with
(non-symmetric) operads. We show that \malg Ss, \mcoalg Ss are
defined as (co)algebras over operads $\Mag^S$. A formula to
compute the dimension series of $\Mag^S$ is given in Proposition
\ref{prop:inversion}.

In Section 4, we introduce \mbialg T Ss where $T\subset
S\subset\NNN$, and we pay special attention to the case $S=T$. The
operads $\magroot T S$ are constructed in Section 5. They are
given (in arity $\geq 2$) by a right operadic ideal (defined
analogously to a right ideal in a ring). We also construct related
forgetful and universal enveloping functors.

In Section 6, we extend results of \cite{B2} to the case of
\mbialg S Ss, and then follow ideas of \cite{LR1} to prove the
main theorem, classifying \mbialg T Ss. In Section 7 we give all
necessary modifications of the previous constructions and the
structure theorem so that they are still valid in the presence of
a unit, provided $S=\{2,\ldots,n\}$ and $T=\{2,\ldots,m\}$.

The structure theorems lead to some pretty combinatorics for the
generating series of the corresponding operads. Some typical
examples are given in Section 8.

\section{\mbialg{T}{S}s}

In this section we define \malg Ss and \mcoalg Ts in order to
construct \mbialg T Ss.
\subsection{\malg{S}s}

\begin{definition}\label{def:S-mag alg}
Let $S$ be a subset of $\NN=\{2,3,\ldots\}$. An \emph{\malg{S}}
$A$ is a vector space endowed with an $n$-ary operation $\mu_n$
for each $n\in S$.
\end{definition}

\begin{definition}\label{def:planar tree}
A \emph{planar tree}  $T$ is a planar graph which is assumed to be
simple (no loops nor multiple edges), connected, rooted and
reduced (no vertices with only one outgoing edge). The set of
planar trees with $n$ leaves is denoted by $\PT_n$, and by $\PT$
we denote the union $\bigcup_{n\geq 1}\PT_n$. The number $n\geq 1$
of leaves will also be called the \emph{degree} of the tree, and
we also allow an empty tree $\emptyset$ of degree 0 later on (in
Section \ref{sectunital}). In low dimensions one gets:
\begin{eqnarray*}
&&\PT_1=\left\{ | \right\} , \ \PT_2=\left\{ \arbreA \right\},
\PT_3= \left\{ \arbreAB \ \arbreBA \arbreAt \right\},\\
&&\PT_4=\left\{ \arbreABC \ \arbreBAC \ \arbreACA \ \arbreCAB
\arbreCBA \right. \\ &&  \left.   \ \arbreABt \ \arbreBAt \
\arbreAAt \ \arbreABttt \ \arbreBAttt \ \arbreAttt \right\} ,
\cdots
\end{eqnarray*}
\end{definition}

The $n$-grafting of $n$ trees is the gluing of the  root of each tree on a new root.
For example the 2-grafting of the two trees $t$ and $s$ is:
$$
\vee_2(t,s):=
\concat \ ,
$$
the 3-grafting of three trees $t$, $s$ and $u$ is:
$$
\vee_3(t,s,u):=
\concatt \ .
$$

\begin{remark}From our definition of a non-empty  planar tree, any $t\in \PT_n$, $n\in \NN$, is of the form
$$ t=\vee_k (t_1,\ldots,t_k)$$ for uniquely determined $k$ and
non-empty trees $t_1,\ldots,t_k$.
\end{remark}

Let $V$ be a vector space. A \emph{labelled tree of degree $n$},
$n\geq 1$,  is a tree $t$ endowed with a labelling of all leaves
by elements $v_1,\ldots, v_n\in V$. We denote such a labelled tree
by $(t,v_1\cdots v_n)$, and represent it as follows:
$$\arbreetiquette $$ One defines the $n$-grafting of labelled
trees by the $n$-grafting of the trees, where one keeps the
labellings on the leaves.

We denote by $\PT^S$ the subset $\PT^S\subset \PT$ that only
contains planar reduced trees where the arities of all their inner
vertices belong to $S\subset \NN$. And by $\PT_n^S$ we denote the
set of planar reduced trees in $\PT^S$ with exactly $n$ leaves.

\begin{example}\label{ex:PT alg}
Let $V$ be a vector space. The vector space $\oplus_{n\geq 1}\mathbb
K[\PT_n^S]\otimes V^{\otimes n}$, spanned by the
labelled planar trees endowed with  the $n$-grafting $\vee_n$ of labelled
trees, for all $n\in S$, is an \malg S.
\end{example}

\begin{definition}\label{def:S-mag coalg}
Let $S\subset \NN$ be a subset. An \emph{\mcoalg{S}} $C$ is a
vector space endowed with an $n$-ary co-operation $\Delta_n$ for
each $n\in S$.
\end{definition}

Given $t\in \PT$ a non-empty tree, it can  be uniquely written as
the $n$-grafting of non-empty trees $t_1,\ldots,t_n$. Then the
\emph{$n$-ungrafting} of planar trees is defined as follows:
$$\wedge_n(t):= t_1\otimes\cdots\otimes t_n$$ and the
$m$-ungrafting $\wedge_m$ of such a tree is zero for $m\neq n$.

\begin{example}\label{ex:PT coalg}
Let $V$ be a vector space. The vector space $\oplus_{n\geq 1
}\mathbb K[\PT_n^S]\otimes V^{\otimes n}$,  spanned by the labelled
planar trees endowed with  the $n$-ungrafting $\wedge_n$ of
labelled trees, for all $n\in S$, is an \mcoalg S.
\end{example}

\section{Recall on operads}
We recall some properties of non-symmetric operads. The operadic
setting is quite useful for our purpose, and it will permit us to
easily derive the combinatorics of the operads involved later on.

\begin{definition}\label{def:Schur functor}
To  a graded vector space
 $$M = (M_0,M_1, \ldots ,M_n, \ldots)\ ,$$ we associate
  its \emph{Schur functor} $\widetilde M : \Vect \rightarrow \Vect$ defined by
\begin{eqnarray*}
\widetilde {M}(V) := \oplus_{n\geq 0} M_n \otimes_\KK V ^{\otimes
n}\ .
\end{eqnarray*}
The Schur functor admits a tensor product and a composition
defined by $\widetilde{M\otimes N}=\widetilde{M}\otimes
\widetilde{N}$ and $\widetilde{M\circ N} =\widetilde{M}\circ
\widetilde{N} $. The tensor product and the composition of Schur
functors can be seen as a Schur functor on  graded vector spaces
denoted by $M\otimes N$ and defined as :
\begin{eqnarray*}
(M\otimes N)_n := \bigoplus_{i+j=n} M(i) \otimes N(j)\ ,
\end{eqnarray*}
and their composite denoted by $M\circ N$ is defined as :
\begin{eqnarray*}
(M\circ N)_n :=\bigoplus_{i_1+\cdots+i_k=n} M_k \otimes(N_{i_1}
\otimes \cdots \otimes N_{i_k}) \ .
\end{eqnarray*}
\end{definition}

\begin{definition}\label{def:nonsymmetric operad}
A \emph{non-symmetric operad} (often abbreviated into operad in
this paper) is a graded vector space $\PP = \{\PP_n\}_{n\geq 0}$
equipped with composition maps $\gamma_{i_1,\cdots,i_n} : \PP_n
\otimes \PP_{i_1} \otimes \cdots \otimes \PP_{i_n} \rightarrow
\PP_{i_1+\cdots+i_n}$ and an element $\Id \in \PP_1$, such that
the transformations of functors $\gamma : \PP  \circ\PP
\rightarrow \PP$ and $\iota : \Id \rightarrow \PP$, deduced from
this data, endows $(\PP,\gamma,\iota)$ with a monoidal structure
on the Schur functor $\PP$.

We will only consider operads $\PP$ with $\PP_0=0$, $\PP_1=\KK\
\Id$ (connected non-symmetric operads).

Let $F$ be the functor from the category of such operads to the
category of graded vector spaces, which forgets the monoidal
structure (and in particular deletes $\PP_1$). Its left adjoint is
the \emph{free non-symmetric operad functor}, $$ M=(0,0,M_2,M_3,
\ldots ,M_n, \ldots)\mapsto \TT(M).$$

With $\PT$ denoting the set of planar reduced trees, the arity
$n$-component $\TT(M)_n$ of the non-symmetric operad $\TT(M)$ can
be identified with the vector space $$\bigoplus_{t\in\PT_n}
M_{i_1}\otimes\cdots\otimes M_{i_k},$$ where $k$ is the number of
internal vertices of $t\in\PT_n$ and $(i_1,\ldots, i_k)$ lists the
arities of the internal vertices.

The operad composition can be described by the plugging (by some
authors also called grafting) of trees onto the leaves (compare
also Definition \ref{def:magS} later).

\end{definition}

\begin{definition}\label{def:operadideal}
Let $\PP$ be a non-symmetric operad (or an operad), and let
$\III(n)$ be a subspace of $\PP_n$ for each $n$.

Then $\III=(\III_n)_{n\in\NN}$ is called a (two-sided) ideal of
$\PP$, if for $$\gamma_{i_1,\cdots,i_n}(p\otimes
q_1\otimes\ldots\otimes q_n)$$ defined in $\PP$, it follows that
$$\gamma_{i_1,\cdots,i_n}(p\otimes q_1\otimes\ldots\otimes q_n) \in\III$$
whenever at least one of the operations $p, q_1,\ldots,q_n$ is in
$\III$.

The quotient of $\PP$ by $\III$ is again a non-symmetric operad
and will be denoted by $\PP\slash\III$.

As in rings, we analogously define right ideals in $\PP$, $\PP$ a
non-symmetric operad. For $\MM=(\MM_n)_{n\in\NN}$ with
$\MM_n\subset \PP_n$, the right ideal $\III:=\MM\circ\PP$
generated by $\MM$ in $\PP$ consists of all elements of the form
$\gamma_{i_1,\cdots,i_n}(p\otimes q_1\otimes\ldots\otimes q_m) $
with $p\in\MM_n$, $q_1,\ldots,q_n\in\PP$ (all $n$).

\end{definition}

\begin{definition}\label{def:algebra over operad}
An \emph{algebra over a non-symmetric operad} $\PP$ (in short:
$\PP$-algebra) is a vector space $A$ equipped with a family of
linear maps $\gamma_{A,n}: \PP_n\otimes A^{\otimes n}\to A$,
called structure maps. In particular, for $V$ a vector space, the
underlying vector space of the free $\PP$-algebra $\PP(V)$ is
associated to $V$ by the Schur functor $\PP$, and the structure
maps $\gamma_{A,n}$ combine to a linear map $\gamma_A : \PP(A)
\rightarrow A$, such that the following diagrams are commutative:
$$ \xymatrix{\PP\circ\PP(A) \ar[r]^{\PP(\gamma_A)}
\ar[d]_{\gamma(\PP(A))}&\PP(A)\ar[d]^{\gamma(A)}\\ \PP(A)
\ar[r]^{\gamma_A} & A \ .} \qquad \xymatrix{A
\ar[r]^{\PP(\iota_A)} \ar[rd]_{=}&\PP(A)\ar[d]^{\gamma(A)}\\
 & A \ .}$$

Morphisms of ${\mathcal P}$-algebras $A\to B$ are linear maps
$\varphi:A\to B$ compatible with the corresponding structure maps
$\gamma_A, \gamma_B$ , i.e.\ such that the following diagram
commutes:
\begin{equation*}
\xymatrix{
          {\mathcal P}(n)\otimes A^{\otimes n}
          \ar[d]_{\Id\otimes\varphi^{\otimes n} }
          \ar[r]^{\ \ \ \ \ \ \ \ \gamma_{A,n}}
          & A
           \ar[d]^{\varphi}\\
            {\mathcal P}(n)\otimes B^{\otimes n}
        \ar[r]_{\ \ \ \ \ \ \ \ \gamma_{B,n}} & B
            }
\end{equation*}

A ${\mathcal P}$-algebra $A$ is called  nilpotent, if for
sufficiently large $n$,
\begin{equation*}
\gamma_A(n)(\mu_n\otimes a_1\otimes\ldots a_{n})=0
\end{equation*}
for all $a_1,\ldots,a_n\in A$, $\mu_n\in\PP_n$.

Analogously defined is the notion of a \emph{coalgebra over a
non-symmetric operad}, it is a vector space $C$ together with a
family of structure maps $\lambda_{C,n}: \PP_n\otimes C\to
C^{\otimes n}$.

Morphisms of ${\PP}$-coalgebras $C\to D$ are $\KK$-linear
maps $\varphi:C\to D$ compatible with the corresponding structure
maps.

Following Loday \cite{L}, we define the notion of connected
coalgebra as follows. A ${\mathcal P}$-coalgebra $C$ is called
connected or (co-)nilpotent, if the following condition holds:
\begin{equation*}
\textrm{ for all } c\in C \textrm{ there exists } r\in \NNN
\textrm{ such that for } n > r, \lambda_{C,n}(\delta_n\otimes
c)=0, \textrm{ all }\delta_n\in\PP_n.
\end{equation*}
For any ${\mathcal P}$-coalgebra $C$, $r\geq 1$, we put
 $$\F {r} C:=\bigcap_{n > r}\ \ \bigcap_{\delta_n\in \PP_n} \left\{x \in C\ \vert \
\delta_n(x) =0\right\}.$$ And $C$ is connected iff $\bigcup_{r
\geq 1}\F {r}C=C$. Note the $\delta_n\in\PP_n$ is any composition
of the generating co-operations such that $\delta_n:C\rightarrow
C^{\otimes n}$.

\end{definition}

The  \malg{S} (Definition \ref{def:S-mag alg}) and  \mcoalg{S}
(Definition \ref{def:S-mag coalg}) can be viewed as algebras over the
non-symmetric operad $\Mag^S$, introduced as follows :

\begin{definition}\label{def:magS}
Let $S\subset \NN$ be a set. Let $M$ be the graded vector space
given by $M_n=\KK$ for $n\in S$ and $M_n=0$ otherwise. Then
$\Mag^S:=\TT(M)$ is the free non-symmetric operad generated by one
$n$-ary generating operation for each $n\in S$.

We can identify the space $\Mag^S_n$ of operations ($n$ inputs)
with the vector space $\KK [\PT^S_n]$. The composition of
operations corresponds to the plugging of trees onto leaves.

Thus we write $\Mag^S(V):=\oplus_{n} \PT_n^S\otimes V^{\otimes n}$
for the (underlying vector space of the) free \malg{S}.
\end{definition}

\begin{remark}
Note that all non-symmetric operads $\Mag^S$ are included in
$\Mag^{\infty}=\Mag^{\NN}$, also denoted by $\Mag^{\omega}$, see
\cite{H}.

Moreover, for every inclusion $T\to S$ there is an inclusion of
operads $\Mag^T\to\Mag^S$.
\end{remark}

\begin{definition}\label{def:operad nilS}

 Let $\Nil^S$ denote the
non-symmetric operad given by the quotient of  $\Mag^S$ with
respect to the (two-sided) ideal consisting of all nontrivial
compositions, i.e.\ any nontrivial composition will be zero in the
quotient.
\end{definition}

The following combinatorial formula is the key to explicit the
dimensions of the spaces $(\Mag^S)_n$ (i.e.\ of the spaces of
operations with $n$ inputs).
\begin{proposition}\label{prop:inversion}
For every $n\in\NN$, the dimension $a_n^S:=\dim ((\Mag^S)_n)$ is
given by the following formula (in the ring of power series in one
variable $x$ over $\KK$):

$$(x - \sum_{i\in S} x^i)\circ (\sum_{n\ge 1} a^S_n x^n)=x,$$
 where
$\circ$ is the composition of power series.

\end{proposition}

\begin{proof}

For a presented quadratic  non-symmetric operad $\PP$, let
$\PP_n^d$ be the space spanned by the $n$-ary operations
constructed out of $d$ generating operations, and let
$f_{\PP}(x,z)$ be defined as follows :
$$f_{\PP}(x,z):= \sum_{{n\ge 1}\atop {d\geq 0}}\dim\PP_n^d \ z^d x^n.
 $$

The operad $\Mag^S$ is free and especially Koszul. In the Koszul
duality theory of operads, as noted for binary quadratic operads
in \cite{GK}, the Poincar\'e series of a Koszul operad and its
Koszul dual are related.

Since we don't work with binary quadratic operads but with general
quadratic operads, we need the following general formula for
Koszul operads, given by B.\ Vallette in \cite{V}(Section 9):
$$f_{\QQQ}(f_{\PP}(x,z),-z)=x, \textrm{ where } \PP^!=\QQQ \textrm{
denotes the Koszul dual operad of } \PP\ .$$

In our case, the Koszul dual of $\Mag^S$ is exactly the operad
$\Nil^S$.

We put $z=1$ and set $f(x):=f_{(\Mag^S)}(x,1)= \sum_{n\ge
1}\dim((\Mag^S)_n) \ x^n.$

We get the following Lagrange inversion type formula: $$ g\circ f
= \Id, \textrm{ for } g(x)=\sum_{{n\ge 1}\atop {d\geq 0}}(-1)^d
\dim((\Nil^S)_n^d)\ x^n. $$ It is clear that the space
$(\Nil^S)_1^{0}=\KK\ \Id$ is one-dimensional. Also it is clear
that
 $\dim
((\Nil^S)_n^{1}) =1$ for all $n\in S$, while $\dim
((\Nil^S)_n^{d}) =0$ for any other $n,d$.

Thus $f(x)$ is the composition inverse of $g(x)=x-\sum_{n\in
S}x^n.$
\end{proof}
We illustrate this formula in some particular cases in section
\ref{sec:combinatorics}.

\section{\mbialg{T}{S}s}
This section introduces the natural notion of \mbialg {T}{S} that
we will classify in section \ref{sec:main thm} under the
connectedness property.

\medskip

The following proposition highlights that the vector space
$\oplus_{n\geq 1}\KK[\PT_n^S]$ admits an algebra structure due to
the $n$-grafting of trees, $\vee_n$ for $n\in S$, and a
coalgebraic structure by endowing it with the $m$-ungrafting of
the trees, $\wedge_m$ for $m\in S$.

\begin{proposition}\label{prop:freecofree mag}
Let $V$ be a vector space.
\begin{enumerate}[i.]
\item
The space $\Mag^S(V)$ endowed with  the $n$-grafting of labelled
trees, $\vee_n$ for all $n\in S$,
 is the free \malg{S}. Explicitly, for $i: V\to\Mag^S(V)$ the inclusion, the
following universal property is fulfilled: \\ Any linear map $f :
V\to A$, where $A$ is any \malg{S}, extends to a unique morphism
of \malg{S}s $\tilde{f} : \Mag^S(V) \to A$ with $f=\tilde{f}\circ
i.$
 $$ \xymatrix{ V \ar[r]^{i} \ar[dr]_{f}&Mag^S(V) \ar[d]^{\tilde f}\\
 & A \ .} $$
\item
The space $\Mag^S(V)$  endowed with  the $n$-ungrafting of
labelled trees, $\wedge_n$ for all $n\in S$,
 is the
cofree \mcoalg{S}. Explicitly, for $p :\Mag^S(V) \to V$ the
projection on $V$, the following universal property is fulfilled:\\
Any linear map $\phi : C\to V$, where $C$ is any connected
\mcoalg{T}, extends to a unique coalgebra morphism $\tilde \phi :
C\to
\Mag^S(V)$: $$ \xymatrix{C \ar[dr]^{\phi} \ar[d]_{\tilde \phi}&\\
\Mag^S(V) \ar[r]^{p} & V \ .} $$
\end{enumerate}
\end{proposition}
Assertions (i) and (ii) of the above Proposition are dual to each
other, and the proofs are by direct inspection, cf.\ similar
proofs in \cite{B2}.

As on the same vector space an algebraic and a coalgebraic
structure are coexisting, a natural question would be to seek for
relations between the operation and co-operations. Note that in
the case of classical bialgebras, the relation (called the Hopf
relation) between the operation and co-operation is determined by
the fact that one is a homomorphism for the other.

\begin{lemma}\label{lemma:rel on MagS}
The relation between the $n$-grafting $\vee_n$ of trees and the
$m$-ungrafting $\wedge_m$ of trees is the following: Let
$t_1,\ldots,t_n\in\PT^S$ be trees, then
\begin{eqnarray*}
\wedge_m\circ\vee_n(t_1,\cdots,t_n)=\left\{
\begin{array}{cl}
0 &\textrm{ if } m\neq n \ ,\\
t_1\otimes\cdots\otimes t_n &\textrm{ otherwise.}
\end{array}
\right.
\end{eqnarray*}
\end{lemma}

\begin{proof}
The ungrafting of trees is defined in example \ref{ex:PT coalg}.
The tree $t=\vee_n(t_1\otimes\cdots\otimes t_n)$ is uniquely
decomposed as product of $n$ trees, by construction of the space
$\KK[\PT^S]$. And by definition the ungrafting depends only on the
arity $n$ of the root. This ends the proof.
\end{proof}
The above proposition motivates the following definition :

\begin{definition}\label{def:T<S mag bialg}
Let $T\subset S\subset \NN$ be two sets. A \emph{\mbialg{T}{S}}
$(\HH,\mu_n,\Delta_n)$ is a vector space $\HH$ endowed with a
\malg{S} and a \mcoalg{T} structure verifying the following
compatibility relation for all $k\in T$ and $l\in S$:
\begin{equation}\label{eq:compatibility relation}
\begin{array}{rcl}
\Delta_k\circ\mu_k =&
\Id\ ,&\\
\Delta_k\circ\mu_l =& 0&\quad \textrm{for } k\neq l\ .
\end{array}
\end{equation}
\end{definition}

\begin{remark}
This is a generalisation of the \mbialg {\m}{\n} structure
introduced in \cite{B2} in the non-unital framework. (The
compatibility relations collapse to the ones above in the
non-unital framework).
\end{remark}

\begin{proposition}\label{prop:structure of malgS on mcoalgT}
There exists a unique family of coproducts $\Delta_m$ with $m\in
S$ such that $\Delta_m(v)=0$ for all $v\in V$ on the free \malg{S}
$\Mag^S(V)$ which makes it into an \mbialg{S}{S} for which $V$ is
primitive. Moreover as a coalgebra $\Mag^S(V)$ is connected.
\end{proposition}

\begin{proof}
The existence is a consequence of lemma \ref{lemma:rel on MagS}.

The uniqueness is due to the compatibility relation. Indeed, let
us construct $\Delta_m:\Mag^S(V)\rightarrow\Mag^S(V)^{\otimes m}$
, for all $m\in S$, as a \mcoalg S co-operation induced by
$\Delta_m(v)=0$ and verifying the compatibility relation
(\ref{eq:compatibility relation}). As any tree $t\in\PT^S$ can be
uniquely viewed as the $n$-grafting of some trees $t_1,\ldots,t_n$
the $m$-ary co-operation on $t$ evaluates to
$\Delta_m(t)=\Delta_m\circ\vee_n(t_1\otimes\cdots\otimes t_n)$.
Therefore by the compatibility relation we get:
\begin{eqnarray*}
\Delta_m\circ\vee_n(t_1,\cdots,t_n)=\left\{
\begin{array}{cl}
0 &\textrm{ if } m\neq n \ ,\\
t_1\otimes\cdots\otimes t_n &\textrm{ otherwise.}
\end{array}
\right.
\end{eqnarray*}
This proves the uniqueness (there is no other choice to construct
the $m$-ary co-operations).

The connectedness is proven in a similar way as in \cite{B2} and
gives the following filtration: $$F_r(\Mag^S(V))=\oplus_{n=1}^{r}
\KK[\PT_n^S]\otimes V^{\otimes n}$$ and $\Prim\Mag^S(V)=V$.
(Remark that we are in the case of a \mbialg S S and not of a
\mbialg T S where $T\neq S$.)
\end{proof}

Remark that the constructed co-operation $\Delta_n$ is exactly
$\wedge_n$ by uniqueness.

\section{The \mrootalg{T}{S}}
In the following, we describe a non-symmetric operad related to
the \mop S $\Mag^S$ and the \mop T $\Mag^T$. It will occur in the
structure theorems.

\begin{definition}\label{def:mag root operad}

Let $\III=(\III_n)_{n\in\NN}$ be the right ideal generated by the
$k$-ary operations $\mu_k$ for all $k\in S\setminus T$ in the
non-symmetric operad $\Mag^S$ (see Definition
\ref{def:operadideal}). We set $(\magroot{T}{S})_1:=\KK\ \Id$ and
$(\magroot{T}{S})_n:=\III_n$ for $n>1$.

\end{definition}

\begin{proposition}

We get a non-symmmetric operad
$$\magroot{T}{S}=((\magroot{T}{S})_n)_{n\in\NN}.$$

Moreover the \mrootop{T}{S} is freely generated by operations
$$\mu_k\circ(\nu_1^S\otimes\cdots\otimes\nu^S_k) \ ,$$ where $k\in
S\setminus T$ and $\nu_i^S$ are operations that are compositions
of generating operations from the \mop S but do not belong to the
\mrootop{T}{S}. The operad structure is given by its composition
map $\gamma:\magroot{T}{S}\circ\magroot{T}{S}\rightarrow
\magroot{T}{S}$ which is the composition of operations  and its
unit map $\iota:\Id\rightarrow\magroot{T}{S} $ which is given by
the operation identity. It is a sub-operad of the \mop S.

\end{proposition}
\begin{proof}
 It is easy to verify the monoidal axioms. Moreover, there
is a unique way of writing the operations as composition of
generating operations (this gives the injectivity of the  map
\mrootop T S to \mop S).
\end{proof}

\begin{definition}\label{def:mag root algebra}
A \emph{\mrootalg{T}{S}} $A$ is an algebra over the free operad
$\magroot{T}{S}$ (i.e. it is a vector space endowed with
operations $\mu_k$, indexed by $S\setminus T$, and
$\mu_k\circ(\nu^S_{1}\otimes\cdots\otimes \nu^S_{k})$ as above).
\end{definition}

We will denote by $\PT^S_{\root\setminus T}$  the set of trees
$t\in \PT^S$ such that the arity of the root is in $S\setminus T$.
And $(\PT^S_{\root\setminus T})_n$ denotes the subset of trees
having exactly $n$ leaves.

We can identify the space $(\magroot T {S})_n$ of operations ($n$
inputs) with the vector space $\KK[(\PT^S_{\root\setminus T})_n]$.
The composition of operations corresponds to the plugging of trees
onto leaves.

Thus we write $\magroot T S(V):=\oplus_{n}
\KK[(\PT^S_{\root\setminus T})_n]\otimes V^{\otimes n}$ for the
underlying vector space of the free \mrootalg T S.

\begin{example}
To illustrate the above identification, we consider the operation
$\mu_2\circ(\Id\otimes\mu_3)\circ(\Id\otimes\Id\otimes\Id\otimes\mu_3)
$, in the space of operations of $\Big(\magroot {\{3,4\}}
{\{2,3,4\}}\Big)_{6}$, which can be depicted as the following tree
: $\arbreCBAt$.
\end{example}

One key part of the proof of the structure theorem is the
construction of the forgetful functor from the algebra structure
of the bialgebra to its primitive part. The universal enveloping
functor is the left adjoint functor to this forgetful functor and
is constructed quite naturally.

Let $T\subset S\subset \NN$ be two sets. Let $A$ be a \malg{S}. By
definition it is equipped with operations $\mu_n$, and
consequently also with all possible compositions of these. By
$$\mu_k\circ({\textrm{composition of operations indexed by S})},$$
where $k$ runs through all elements of $S\setminus T$, we
construct an infinite family of operations (contained in the
family of all operations).

\begin{proposition}\label{prop:forgetful functor}
The above construction defines a functor
\begin{eqnarray*}
()_{\magroot{T}{S}}:&\{\textrm{\malg{S}}\}\rightarrow\{\textrm{\mrootalg{T}{S}}\}\\
&(A,\{\mu_k\}_{k\in S})\mapsto\\
&(A,\{\mu_k\circ({\textrm{composition of operations indexed by
S}})\}_{k\in S\setminus T})\ ,
\end{eqnarray*}
namely the forgetful functor.
\end{proposition}
\begin{proof}
Definition \ref{def:mag root algebra} ensures that the constructed
algebra $$(A,\{\mu_k\circ({\textrm{composition of operations
indexed by S}})\}_{k\in S\setminus T})$$ is a \mrootalg{T}{S}.
\end{proof}

We construct a functor
$$U_{\magroot{T}{S}}:\{\textrm{\mrootalg{T}{S}}\}\rightarrow\{\textrm{\malg{S}}\}
$$ as follows: For any \mrootalg{T}{S}
$(A,\{\tilde\mu_k\circ(\tilde\nu_1^S\otimes\cdots\otimes\tilde\nu_k^S)_{k\in
S\setminus T}\})$, let the \malg{S} $U_{\magroot{T}{S}}(A)$ be
constructed as the free \malg{S} over the vector space $A$
quotiented by the relations which identify
$\tilde\mu_k\circ(\tilde\nu_1^S\otimes\cdots\otimes\tilde\nu_k^S)$,
for $k\in S\setminus T$, with the corresponding operations $
\mu_k\circ(\nu_1^S\otimes\cdots\otimes\nu_k^S)$ in the free
\malg{S}.

\begin{proposition}\label{prop:adjoint functor}
The universal enveloping functor
$$U_{\magroot{T}{S}}:\{\textrm{\mrootalg{T}{S}}\}\rightarrow\{\textrm{\malg{S}}\} $$
is left adjoint to the forgetful functor
$$()_{\magroot{T}{S}}:\{\textrm{\malg{S}}\}\rightarrow\{\textrm{\mrootalg{T}{S}}\}\ .$$
\end{proposition}

\begin{proof}
Let $A$ be a \malg{S} and $B$ be a \mrootalg{T}{S}. On the one
hand, let $f:B\rightarrow (A)_{\magroot{T}{S}}$ be a
\mrootalg{T}{S} morphism. It determines uniquely  a morphism of
\malg{S} $\Mag^S(B)\rightarrow A$  since $A$ is a \malg{S} and
$\Mag^S$ is the free \malg S.

On the other hand, let $g:U_{\magroot{T}{S}}(B)\rightarrow A$ be a
\malg S morphism. The construction of the universal enveloping
functor $U_{\magroot{T}{S}}$ gives that the map $B\rightarrow
U_{\magroot{T}{S}}(B)$ is a \mrootalg{T}{S} morphism. Hence the
composition with $g$ gives a \mrootalg{T}{S} morphism
$B\rightarrow (A)_{\magroot{T}{S}}$.
\end{proof}

\begin{corollary}\label{cor:universal enveloping functor of the free mrootalgTS iso to free malgS}
The universal enveloping functor of the free \mrootalg {T}{S} is isomorphic to the free \malg S :
$$U_{\magroot{T}{S}}(\magroot{T}{S}(V))\cong \Mag^S(V)$$
\end{corollary}

\begin{proof}
Note that the underlying vector space is preserved under the
forgetful functor $()_{\magroot{T}{S}}$. As the universal
enveloping functor $U_{\magroot{T}{S}}$ is left adjoint to the
functor $()_{\magroot{T}{S}}$ and that the left adjoint of the
forgetful functor $$\{\textrm{\mrootalg T S}\}\rightarrow\Vect$$
is the free functor their composite is left adjoint to the
forgetful functor $$\{\textrm{\malg S}\}\rightarrow\Vect.$$
Therefore it is the free functor $\Vect\rightarrow\{\textrm{\malg
S}\}$. $$ \xymatrix{\{\textrm{\malg S} \}\ar[r]
&\ar@/^{0.5cm}/^{U}[l]\{\textrm{\mrootalg T S}\}\ar[r]&\Vect
\ar@/^{0.5cm}/^{free}[l]\ar@/_{1cm}/_{free}[ll]} $$
\end{proof}

\begin{corollary}\label{cor:U(V) is a malgS}
For any \mrootalg{T}{S} $A$, $U_{\magroot{T}{S}}(A)$ is a connected \mbialg{S}{S}.
\end{corollary}

\begin{proof}
It is a consequence of the above Corollary \ref{cor:universal
enveloping functor of the free mrootalgTS iso to free malgS} and
Proposition \ref{prop:structure of malgS on mcoalgT}.
\end{proof}

\section{Main theorem}\label{sec:main thm}
In this section, we state the classification theorem for connected
\mbialg T S named the structure theorem. In order to prove it, we
focus on the construction of the algebra of the primitive
elements, and on the rigidity theorem in the particular case of an
\mbialg SS.

\begin{definition}\label{def:primitive element}
Let $C$ be a \mcoalg{T}. The vector space of \emph{primitive
elements} $\Prim C$ is defined as
$$\Prim C:= \cap_{n\in T} \left\{x \in \mathcal{H}\ \vert \   \Delta_n(x)=0 \right\}$$
\end{definition}

The following theorem gives the classification of the connected
\mbialg{T}{S}s.
\begin{theorem}\label{thm:structure thm}
Let $T\subset S\subset \NN$ be two sets. If $\HH$ is a
\mbialg{T}{S} over a field $\KK$, then the following are
equivalent:
\begin{enumerate}[i.]
\item \label{thm:enum:connectedness} $\HH$ is a connected \mbialg{T}{S}.
\item \label{thm:enum:isomorphism} $\HH$ is isomorphic to $U_{\magroot{T}{S}}(\Prim \HH)$ as a \mbialg{T}{S}.
\item \label{thm:enum:cofree}$\HH$ is cofree among the connected \mcoalg{T}.
\end{enumerate}
\end{theorem}

To prove the above theorem, we have to make the primitive part of
a \mcoalg T explicit.
\begin{proposition}\label{prop:primitive part}
The primitive part of  a  \mbialg{T}{S} is generated by the
operations $\mu_k$ and all its composites
$\mu_k\circ(\nu_{1}\otimes\cdots\otimes\nu_{k})$  for $k\in
S\setminus T$, and $\nu_i$ are operations generated by operations
indexed by $S$.
\end{proposition}
\begin{proof} By definition the compatibility relation
(\ref{eq:compatibility relation}) gives $\Delta_k\circ \mu_l=0$
for any $k\neq l$. It comes naturally that $\mu_k$, for $k\in
S\setminus T$, is primitive. Indeed, first we have
$\Delta_m\circ\mu_k=0$ for all $m\in T$ and therefore any
composite of the generating co-operations $\delta$  composed with
$\mu_k$ will be zero. This gives the primitiveness property. And
this is also true for all its composites
$\mu_k\circ(\mu_{i_1}\otimes\cdots\otimes\mu_{i_k})$ such that
$k\in S\setminus T$ and $\mu_{i_1},\ldots,\mu_{i_k}$ composition
of operations indexed by $S$.

As the \mop S is free, any operation can be uniquely written as
composition of the generating operations. Therefore there is a
unique way to write an operation as
$\mu_k\circ(\mu_{i_1}\otimes\cdots\otimes\mu_{i_k})$ where $\mu_k$
is a generating operation and $\mu_{i_j}$ are operations of the
\mop S. Suppose that $k\in T$. Then the $k$-ary co-operation
$\Delta_k$ composed with
$\mu_k\circ(\mu_{i_1}\otimes\cdots\otimes\mu_{i_k})$ is not zero
but equals $\mu_{i_1}\otimes\cdots\otimes\mu_{i_k}$. And
$\mu_k\circ(\mu_{i_1}\otimes\cdots\otimes\mu_{i_k})$ is not
primitive in this case.

Therefore the only primitive operations are the operations $\mu_k$
and all its composites
$\mu_k\circ(\mu_{i_1}\otimes\cdots\otimes\mu_{i_k})$ such that
$k\in S\setminus T$ and $i_1,\ldots,i_k\in S$, which ends the
proof.
\end{proof}

\begin{example} We consider  \mbialg {\{2,4\}}
{\{2,3,4,5\}}s, and especially the free algebra. Thus the
 tree $t:=\arbreABt$ represents an element which is primitive. Indeed, the following
holds
\begin{eqnarray*}
\delta_2(t)=\delta_4(t)=0\ ,
\end{eqnarray*}
by the compatibility relations (\ref{eq:compatibility relation}).
\end{example}

\begin{proposition}\label{prop:structure mag root on primitive part}
Let $(\HH,\mu_n,\Delta_m)$ be a \mbialg{T}{S}. Its primitive part
admits a \mrootalg{T}{S} structure.
\end{proposition}

\begin{proof}
The generating operations of the primitive part verify the
conditions of definition \ref{def:mag root algebra}.
\end{proof}

\begin{corollary}\label{cor:primitive part of free malg}
The primitive part of the free \malg S $\Mag^S(V)$ over a  vector space $V$ is exactly $\magroot T S(V)$.
\end{corollary}

\begin{proof}
The first remark is that any tree $t$ which admits a root of arity
$n\in S\setminus T$ is primitive. This
 is due to the fact that the tree can be uniquely rewritten as
 $\mu_n(t_1\otimes\cdots\otimes t_n)$ for certain non-empty trees $t_1,\ldots, t_n$.
 And  by the compatibility relation gives the relation $\wedge_m\circ\vee_n=0$
 (Lemma \ref{lemma:rel on MagS}). This prove that these trees are primitive elements.
Similarly, trees such that their root is of arity $n\in T$ are not
primitive as there exists $\wedge_n$ such that
$\wedge\circ\mu_n=\Id$. By the isomorphism identifying the trees
with root of arity $n\in S\setminus T$ with \mrootalg T S the
proof is completed (Definition \ref{def:mag root algebra}).
\end{proof}

Then, we focus on the particular case of the structure theorem
where the sets $T$ and $S$ are equal. The rigidity theorem is a
classification of connected \mbialg S S where the primitive part
admits a structure of vector space. The results of  \cite{B2}
where $S=\{2,\cdots, n\}$ can be extended to the context of
\mbialg S S with very few changes.

\begin{definition}\label{def:completed algebra}
The \emph{completed \malg S}, denoted by $\Mag^S(\KK)^{\wedge}\ ,$
is defined by $$\Mag^S(\KK)^{\wedge}~\!\!=~\!\!\prod_{n\geq
0}Mag_n^S \ .$$  This definition allows us to define
 formal power series of trees in $Mag^S(\KK)^{\wedge}$, i.e.\ we consider formal power series
  in the non-associative variable $\vert$,
  where $\vert$ denotes the generator of the one-dimensional space $\KK=Mag_1^S(\KK)$, cf. \cite{G}.
  \end{definition}

\begin{definition}\label{def:convolution}
Let $\HH$ be an  \mbialg S S, $n\in S$, and let $f_1,\ldots ,f_n$
be linear maps $\HH\to \HH$.
  The \emph{$n$-convolution} of $f_1,\ldots ,f_n$ is the
linear map defined by:
\begin{equation*}
 \star_n (f_1\cdots f_n) := \mu_n \circ (f_1 \otimes \cdots \otimes f_n) \circ \Delta_n  \ .
 \end{equation*}

We define a map $\chi$ from the set of trees to vector space of
operations of the \mbialg S S: to the $n$-th corolla
$t_n=\vee(\underbrace{|,\ldots,|}_{n})$ we associate the operation
$\star_n$. As any other tree can be seen as a grafting of corollas
of degree $k$, it is associated with the composition of the
respective operations $\star_k$. We denote $\star_t:=\chi(t)$.
\end{definition}

\begin{theorem}\label{thm:S=T rigidity thm}
Any connected \mbialg S S $\HH$ is isomorphic to
$$\Mag^S(\Prim\HH):=(\Mag^S(\Prim\HH),\vee_k,\wedge_k)_{k\in
S}$$ where $\vee_k$ is $k$-grafting and $\wedge_k$ is the
$k$-ungrafting.
\end{theorem}

\begin{proof}

Since the proof is similar to the proof of the analogous result
for $m$-magmatic bialgebras in \cite{B2}, we give only a sketch of
the proof.

Suppose that $\HH$ is connected we prove the isomorphism $\HH
\cong \Mag^S(\Prim \HH)$ by explicitly giving the two inverse
maps. We define the \mcoalg S morphism: $G: \HH \rightarrow
\Mag^S(\Prim\HH)$ as the unique extension of the following linear
map :
\begin{equation*}
x\mapsto x-\sum_{n\in S}\star_n\circ \Id^{\otimes n}(x) \ .
\end{equation*}
Note that $e:\HH \rightarrow \HH$ defined as
\begin{equation*}
e:=\Id-\sum_{n\in S}\star_n\circ \Id^{\otimes n}
\end{equation*}
plays the role of a projector on the primitive part.

Then we define the \malg S morphism $F: \Mag^S(\Prim\HH)
\rightarrow \HH$  as the unique extension of the linear map:
\begin{equation*}
x\mapsto \sum {\star_ t}(\Id)(x) \ ,
\end{equation*}
where the sum extends on all non-empty planar trees
$t\in\oplus_{n\geq 1}\KK[\PT_n^S]$.

Moreover, denote by $y$ the generator of $Mag^S(\mathbb K)$,
$y:=|$, and by $y^n := \vee_n\circ y^{\otimes n}$. We define
$g(y):=y -\sum_{n\in S}y^n,$ and $f(y):= \sum t $, where the sum
extends on all planar trees $t\in\bigcup_{n\geq
 1} \PT_n^S$. Using that these two preceding maps are inverse
with respect to the composition of such power series one can show
that:
\begin{eqnarray*}
F \circ G =\Id_\HH\ ,\quad
G \circ F =\Id_{\Mag^S(\Prim \HH)}\ .
\end{eqnarray*}

 Remark that $\HH \cong Mag^S(\Prim \HH)$ is a \mbialg S S.
Indeed, we have the two following properties :
\begin{eqnarray*}
\vee_n(G(x_1),\ldots,G(x_n))&=&G\circ F(\vee_n(G(x_1),\ldots,G(x_n)))\\
&=&G\circ\mu_n(F\circ G(x_1),\ldots,F\circ G(x_n))\\
&=&G\circ\mu_n(x_1,\ldots,x_n)\\
\Delta_n(F(x))&=&((F\circ G)\otimes\cdots\otimes (F\circ G))\circ\Delta_n(F(x))\\
&=&(F\otimes\cdots\otimes F)\circ\wedge_n(G\circ F(x))\\
&=&(F\otimes\cdots\otimes F)\circ\wedge_n(x) \ ,
\end{eqnarray*}
which proves that $F$  is moreover \mcoalg S morphism (resp. $G$
is a \malg S morphism and hence an \mbialg S S morphism.

\end{proof}

Now we can  prove the structure theorem stated at the
beginning of the section.
\medskip

Before starting the proof we recall that the primitive space only
depends on the co-operations (there is no unit here), i.e. on the
coalgebraic structure.
\begin{proof}[Proof of theorem \ref{thm:structure thm}]
We prove the following implications
\ref{thm:enum:connectedness}$\Rightarrow$\ref{thm:enum:cofree}$\Rightarrow$\ref{thm:enum:isomorphism}$\Rightarrow$\ref{thm:enum:connectedness}.
\begin{enumerate}[ ]
\item \ref{thm:enum:connectedness}$\Rightarrow$\ref{thm:enum:cofree}.
If $\HH$ is a \mbialg{T}{S} then by theorem \ref{thm:S=T rigidity
thm} $\HH$ is isomorphic to $\Mag^T(\Prim \HH)$ as a
\mbialg{T}{T}. Therefore $\HH$ is cofree.
\item \ref{thm:enum:cofree}$\Rightarrow$\ref{thm:enum:isomorphism}.
If $\HH$ is cofree, then it is isomorphic to $\Mag^T(\Prim \HH)$
and $\Prim \HH$ is a \mrootalg{T}{S}. Moreover $\Mag^T(\Prim\HH)$
is a \malg{S} by endowing the vector space $\Mag^T(\Prim\HH)$ with
the products $\mu_n$ inherited by the products on $\HH$. By
adjunction (proposition \ref{prop:adjoint functor}), the inclusion
map $\Prim \HH\rightarrow \Mag^T(\Prim \HH)$ gives rise to a
\malg{S} morphism: $$\Phi: U_{\magroot{T}{S}}(\Prim
\HH)\rightarrow \Mag^T(\Prim \HH)\ .$$ On the other hand, the
inclusion $\Prim\HH\rightarrow U_{\magroot{T}{S}}(\Prim \HH)$
admits a unique extension $$\Psi:\Mag^T(\Prim \HH)\rightarrow
U_{\magroot{T}{S}}(\Prim \HH)\ .$$ Check that both maps are
inverse one to the other.
\item \ref{thm:enum:isomorphism}$\Rightarrow$\ref{thm:enum:connectedness}.
This is corollary \ref{cor:U(V) is a malgS}.
\end{enumerate}
\end{proof}

\begin{corollary}\label{corollary:pbw}
The following holds for the particular \mbialg T S
$(\Mag^S(V),\vee_n,\wedge_m)_{n\in S,m\in T}$:
$$\Mag^S(V)=\Mag^T\circ\magroot T S(V) \ .$$
\end{corollary}
\begin{proof}
The free \malg S $\Mag^S (V)$ over a vector space $V$ is connected
as a \mcoalg S (proposition \ref{prop:structure of malgS on
mcoalgT}). Similarly it can be proven that it is connected as a
\mcoalg T. Indeed, the filtration induced by the connectedness is:

 $$F_r \magroot T S (V)=
\bigoplus_{l=1}^{r}\KK[(\PT^S_{\root\setminus T})_l]\otimes
V^{\otimes l}$$
 The proof is done by
contradiction and descending induction as in \cite{B2}.

Then applying the structure theorem \ref{thm:structure thm}, we get that
$$\Mag^S(V)=\Mag^T(\Prim \Mag^S(V)) \ .$$
Corollary \ref{cor:primitive part of free malg} ends the proof.
\end{proof}

\section{The unital version of the $m,n$-magmatic
bialgebra}\label{sectunital}
 In this section, we treat a
particular case which can be extended to the unital framework:
$S=\{2,\ldots,n\},\ T=\{2,\ldots,  m\}$. So first we have to
redefine the notion of \malg {\n}, \mcoalg \m, \mbialg {\m} {\n}
in the unital framework. In this case the compatibility relation
between the operations and co-operations is different. But most of
the conclusions of the non-unital framework are preserved.

\begin{definition}\label{def: m-malg}
A \emph{unital \malg \m} $A$ is a   vector space endowed with
$k$-ary operations $\mu_k$, one for each $k$ with $2\leq k\leq m$,
 which are unitary, that is~:

there exists a unit map $\eta:\mathbb K\rightarrow A$, $a\mapsto
a\cdot 1$, such that
 $$
\mu_k(x_1,\cdots,x_k)=\mu_{k-1}(x_1,\cdots,x_{i-1},x_{i+1},\cdots,x_k)
\quad\textrm{ where } x_i=1 \textrm{ and }  x_j\in A, \ \forall j
\ . $$ Diagrammatically this condition is the commutativity of: $$
\produit $$ where $\eta:\mathbb K\rightarrow A$ is the unit map.

From now on we denote by $\Mag^{\m}(V)$ the free unital \malg \m\
on $V$.

\end{definition}

\begin{definition}\label{def:m-mcoalg}

A vector space $C=\KK 1\oplus\overline{C}$ is called a
 \emph{unital \mcoalg \m} $C$ if it is
endowed
 with $k$-ary co-operations $\Delta_k:C\to C^{\otimes k}$, one for each $2\leq k\leq
 m$,
which are co-unitary, i.e. the following diagram is commutative:
$$ \coproduit $$
 where  $\epsilon:C\longrightarrow \mathbb K$ is
the augmentation (and $\Delta_1:=\Id$).
\end{definition}

\begin{definition}\label{def:reduced coproduct}
Let $C=\KK 1\oplus\overline{C}$ be a unital \mcoalg{\m}. The set
of $(p,q)$-shuffles is denoted $\Sh(p,q)$. It is the set of
permutations of $(1,\ldots,p,p+1,\ldots p+q)$ such that the image
of the elements $1$ to $p$ and of the elements $p+1$ to $p+q$ are
in increasing order. The reduced co-operations $\delta_k$ for
$k\in \m$ are defined as follows :
\begin{eqnarray*}
&&\delta_1(x):=x\\ &&\delta_k(x):=\Delta_k(x)- \sum_{l\leq k-1}
\sum_{\sigma\in \Sh(l,k-l)}\sigma\circ (\delta_l(x),1^{\otimes
k-l})\ .
\end{eqnarray*}
\end{definition}

\begin{definition}\label{def:unital primitive element}
Let $C$ be a  unital \mcoalg{\m}. The vector space of
\emph{primitive elements} $\Prim C$ is defined as
\begin{eqnarray*}
&\Prim C :=\\ &\cap_{1< k\leq m} \left\{x \in C\ \vert \
\delta_k(x)=0 \textrm{ for any } k \textrm{-ary reduced
co-operation } \delta_k:C\to C^{\otimes k}
\right\}.\end{eqnarray*} Moreover $C$ is said to be connected if
it verifies: $$C=\cup_r \F r C \ ,$$ with
 $$\F {0} C=\KK,\ \F {r} C:=\bigcap_{k > r} \left\{x \in C\ \vert \
\delta_k(x) =0 \textrm{ for any } k \textrm{-ary reduced
co-operation} \right\}.$$
\end{definition}

\begin{example}\label{exp:unitalbialg}(cf.\ \cite{B2})
The vector space of planar trees $\oplus_{l\geq 0}
\KK[(\PT^{\n})_l]$, where $\PT_0:=\left\{ \emptyset \right\}$, can
be endowed with an \malg {\n} structure and an \mcoalg {\n}
structure with the grafting and the ungrafting of trees defined as
follows : $$\wedge_k(t):=\sum t_1\otimes\cdots\otimes t_k$$ where
the sum extends on all the ways to write $t$ as
$\vee_k(t_1,\ldots,t_k)$, where $t_i\in\bigcup_{l\geq 0}\PT_l$.
 This can be made
explicit as  follows, as for the tree $t$ there is a unique way to
be written as a grafting $t=\vee_k(t_1,\ldots,t_k)$, where
$t_i\neq \emptyset$ for all $i$:
\begin{eqnarray*}
&&\wedge_k(t):= \left(\deconcatk
\right)+\sum_{i=0}^{k-1}\emptyset^{\otimes i}\otimes t\otimes
\emptyset^{\otimes k-i-1}\ ,\\ &&\wedge_l(t):=\left\{
\begin{array}{ll} \sum_{i=0}^{l-1}\emptyset^{\otimes i}\otimes
t\otimes \emptyset^{\otimes l-i-1}\ , & \textrm{if $l<k$}\\
\\
\sum_{i=0}^{l-1}\emptyset^{\otimes i}\otimes t\otimes
\emptyset^{\otimes l-i-1}+\\
\sum_{i_1+\cdots+i_{k+1}=l-k}\emptyset^{\otimes i_1}\otimes t_1
\otimes \emptyset^{\otimes i_2}\otimes \cdots\otimes t_k\otimes
\emptyset^{\otimes i_{k+1}}, & \textrm{if $l>k$}\end{array}
\right.\\ &&\wedge_l(|):=\sum_{i=0}^{l-1} \emptyset^{\otimes
i}\otimes |\otimes \emptyset^{l-i-1},\\
&&\wedge_l(\emptyset):=\emptyset^{\otimes l} \ .
\end{eqnarray*}
\end{example}

\begin{definition}\label{def:mbialg m n}
A \emph{unital \mbialg \m \n} $(\mathcal{H},\mu_k,\Delta_l)$,
where $2\leq k\leq m$, and $2\leq l\leq m$ is a vector space
$\mathcal{H}=\overline{ \mathcal{H}} \oplus \mathbb{K} 1$ such
that:\\ 1) $\mathcal{H}$ admits a $\n$-magmatic algebra structure
with $l$-ary operations denoted~$\mu_l$\!\!~,\\ 2) $\mathcal{H}$
admits a $\m$-magmatic coalgebra structure with $k$-ary
co-operations denoted~$\Delta_k$,\\ 3) $\mathcal{H}$ satisfies the
following ``compatibility relation'':

\begin{equation}\label{eq:unitaly compatibility relation}
\begin{array}{l}
\Delta_l\circ\mu_l (x_1\otimes\cdots\otimes x_l)=
x_1\otimes\cdots\otimes x_l+\sum_{i=0}^{l-1}1^{\otimes i}\otimes
\underline x \otimes 1^{\otimes l-i-1}\ ,\\ \Delta_k\circ\mu_l
(x_1\otimes\cdots\otimes x_l)=\\ \qquad\qquad\qquad\left\{
\begin{array}{ll} \sum_{i=0}^{k-1}1^{\otimes i}\otimes \underline
x \otimes 1^{\otimes k-i-1}\ , & \textrm{if $k<l$}\\
\\
\sum_{i=0}^{k-1}1^{\otimes i}\otimes \underline x \otimes
1^{\otimes k-i-1}+\\ \sum_{i_1+\cdots+i_n+1=k-l}1^{\otimes
i_1}\otimes x_1 \otimes 1^{\otimes i_2}\otimes \cdots\otimes
x_l\otimes 1^{\otimes i_{l+1}} & \textrm{if $k>l$}\end{array}
\right.
\end{array}
\end{equation}
$\forall  \underline x :=\mu_l ( x_1\otimes\cdots\otimes x_l) \
and \  x_1,\ldots,x_l \in \overline{\mathcal{H}}$ and $2\leq k\leq
m$, $2\leq l\leq n$.
\end{definition}

\begin{remark}\label{rem:compatibility relation unified}
The compatibility relation (\ref{eq:unitaly compatibility
relation}) can be unified as follows:

for any elements $x_1,\ldots,x_m\in \HH,$
$$\Delta_k\circ\mu_l(x_1\otimes\cdots\otimes x_l)=\sum_{\mu_k(y_1\otimes\cdots\otimes
y_k)=\mu_l(x_1\otimes\cdots\otimes x_l)}y_1\otimes\cdots\otimes
y_k \ .$$
\end{remark}

\begin{proposition}\label{prop:structure of malg n on mcoalg m}
There exists a unique family of coproducts $\Delta_m$, $m\in \n,$
such that $$\Delta_m(1)=1^{\otimes m},\
\Delta_m(v)=\sum_{i=0}^{m-1} 1^{\otimes i}\otimes v\otimes
1^{\otimes m-i-1} \textrm{ for all } v\in V,$$
 on the free unital
\malg{\n} on $V$ which makes it into a connected unital
\mbialg{\n}{\n} for which $V$ is primitive.
\end{proposition}

\begin{proof}
The existence is due to example \ref{exp:unitalbialg}, where we
identify the unit $1$ with the empty tree $\emptyset$. The
uniqueness is due to the compatibility relation. The proof is
similar to the non-unital framework. Indeed, let us construct
$\Delta_m:\Mag^{\n}(V)\rightarrow\Mag^{\n}(V)^{\otimes m}$ , for
all $m\in \n$ as a \mcoalg \n\ co-operation induced by :
\begin{eqnarray*}\
&&\Delta_m(1)=1^{\otimes m} \ ,\\
&&\Delta_m(v)=\sum_{i=0}^{m-1} 1^{\otimes i}\otimes v\otimes
1^{\otimes m-i-1} \ ,
\end{eqnarray*}
and verifying the compatibility relation (\ref{eq:unitaly
compatibility relation}). It suffices to consider $V=\KK$. As any
tree $t\in\PT^{\n}\setminus\PT^{\n}_1$ can be uniquely viewed as
the $k$-grafting of trees $t_1,\ldots,t_k$ with $t_i\neq
\emptyset$ for all $i$ (for some $k>1$), the $m$-ary co-operation
on $t$ evaluates to
$\Delta_m(t)=\Delta_m\circ\vee_k(t_1\otimes\cdots\otimes t_k)$.
Therefore the compatibility relation induces the value of
$\Delta_m(t)$. This proves the uniqueness (there is no other
choice to construct the $m$-ary co-operations). Note that the
counit relation is verified. Indeed, let $t_1,\ldots,
t_k\in\PT^{\n}$ non empty, then:
\begin{eqnarray*}
&&\wedge_m\circ\vee_{k+1}(t_1\otimes\cdots\otimes\underbrace{\emptyset}_{(i)}\otimes\cdots\otimes
t_k)=\\ &&=\sum_{\vee_{m}(s_1\otimes\cdots\otimes
s_m)=\vee_{k+1}(t_1\otimes\cdots\otimes\emptyset\otimes\cdots\otimes
t_k)}s_1\otimes\cdots\otimes s_m\\
&&=\sum_{\vee_{m}(s_1\otimes\cdots\otimes
s_m)=\vee_{k}(t_1\otimes\cdots\otimes t_k)}s_1\otimes\cdots\otimes s_m\\
&&=\wedge_m\circ\vee_{k}(t_1\otimes\cdots\otimes t_k).
\end{eqnarray*}
This proves that the counit relation is verified.

We get $\Prim\Mag^{\n}(V)=V$ and the connectedness, with the
following filtration: $$F_r(\Mag^{\n}(V))=\KK
1\oplus\bigoplus_{m=0}^{r} \mathbb K[\PT^{\n}_m]\otimes V^{\otimes
m}.$$
\end{proof}

\begin{lemma}\label{lem:mrootalg m n}
Every \mrootalg{\m}{\n} $\bar A$ admits a unit 1, such that
$A:=\KK1\oplus\bar A$ is equipped with $\magroot {\m}
{\n}$-operations.
\end{lemma}
\begin{proof}
In the particular case where $(T;S)=(\m,\n)$, insertion of units
into a $k$-ary operation from $\magroot{\m}{\n}$ yields $k'$-ary
operations from $\magroot{\m}{\n}$ with $k'<k$.
\end{proof}

\begin{proposition}\label{prop:unital forgetful functor}
There is a forgetful functor (in the unital case) given by
\begin{eqnarray*}
()_{\magroot{\m}{\n}}:&\{\textrm{\malg{\n}}\}\rightarrow\{\textrm{\mrootalg{\m}{\n}}\}\\
&(A,\{\mu_k\}_{k\in \{2,\cdots,n\}}) \mapsto\\ &
(A,\{\mu_k\circ(\textrm{operation from the \malg \n})\}_{k\in
\{m+1,\ldots, n\}}).\\
\end{eqnarray*}
\end{proposition}

\begin{proof} The proof in the non-unital framework has to be completed
by the verification of the unital relations. Indeed, the
generating operations $$\{\mu_k\circ(\textrm{operation from the
\malg \n})\}_{k\in \{m+1,\ldots, n\}}$$ are again given by a right
ideal of an \mop \n.\end{proof}

Next we focus on the construction of the universal enveloping
functor in this context. Now, for $A$ a \mrootalg{\m}{\n}, we
construct \malg{\n} $U_{\magroot{\m}{\n}}(A)$ as the quotient of
the free \malg{\n} $\Mag^n(A)$ by the relations which identify
$\tilde\mu_k\circ(\ldots)$ with the corresponding operations in
the free \malg{\n}, and furthermore by the relations identifying
the unitary relations $\widetilde{\RRR}$ verified in $A$ with the
unitary relations ${\RRR}$ verified in $\Mag^n(A)$.

\begin{proposition}\label{prop:unital adjoint functor}
The universal enveloping functor
$$U_{\magroot{\m}{\n}}:\{\textrm{\mrootalg{\m}{\n}}\}\rightarrow\{\textrm{\malg{\n}}\} $$
is left adjoint to the forgetful functor
$$()_{\magroot{\m}{\n}}:\{\textrm{\malg{\n}}\}\rightarrow\{\textrm{\mrootalg{\m}{\n}}\}\
.$$
\end{proposition}

\begin{proof}
The proof is similar to the proof of Proposition \ref{prop:adjoint
functor} except that we must verify that the morphism
$f:B\rightarrow (A)_{\magroot{\m}{\n}}$ passes through the quotient.
Indeed  passing it to the quotient  gives the \malg {\n} morphism
$U_{\magroot{\m}{\n}}(B)\rightarrow A$ as the image
$\RRR(x_1,\ldots,x_n)$ and $\widetilde\RRR(x_1,\ldots,x_n)$ are
the same, namely $\RRR(f(x_1),\ldots,f(x_n))$ and
$\widetilde\RRR(f(x_1),\ldots,f(x_n))$.
\end{proof}

\begin{corollary}\label{cor:universal enveloping functor of the free mrootalg m n iso to free malg n}
The universal enveloping functor of the free unital \mrootalg
{\m}{\n} is isomorphic to the free unital \malg {\n} :
$$U_{\magroot{\m}{\n}}(\magroot{m}{n}(V))\cong \Mag^n(V)$$
\end{corollary}

\begin{theorem}\label{thm:unital structure thm}
Let $m,n\in \NN$ with $m\leq n$. If $\HH$ is a unital \mbialg {\m}
{\n} over a field $\KK$ of characteristic zero, then the following
are equivalent:
\begin{enumerate}[i.]
\item \label{thm:enum:unital connectedness} $\HH$ is a connected unital
\mbialg \m \n. \item \label{thm:enum:unital isomorphism} $\HH$ is
isomorphic to $U_{\magroot{\m}{ \n}}(\Prim \HH)$ as a unital
\mbialg \m \n. \item \label{thm:enum:unital cofree}$\HH$ is cofree
among the connected unital \mcoalg \m s.
\end{enumerate}
\end{theorem}
The proof is similar to the non-unital case treated above.
\section*{}
\begin{remark}

Unital \malg S, \mcoalg S, \mbialg T S would be better described
by pseudo-operads with a zero-ary operation and without
associativity of composition for the cases where this operation is
involved. The exploration of this generalized operadic setting is
still open, compare Borisov-Manin \cite{BM}.

\end{remark}

\section{Combinatorics of the operads $\magroot{T}{S}$}\label{sec:combinatorics}

\begin{lemma}\label{lem:combinatorics}
The generating series $f_{\magroot{T}{S}}(x)= \sum_{n\ge
1}\dim\magroot{T}{S} x^n$ of the non-symmetric operad
$\magroot{T}{S}$ is given by $$x+\sum_{k\in S\setminus
T}(f_{\Mag^S}(x))^k,$$
 where $f_{\Mag^S}(x)$ is the generating
series of the non-symmetric operad $\Mag^S$.
\end{lemma}

\begin{proof}
The dimension of the space $(\magroot{T}{S})_n$ of operations with
$n$ arguments is given by the number of planar trees with $n$
leaves, which fulfill the property, that the arity $k$ of the root
is an element of $S\setminus T$ and the arity of each other
internal vertex is an element of $S$. Equivalently, every tree
representing an operation in $\magroot{T}{S}$ is either the tree
which consists only of the root, or is given by the choice of
$k\in S\setminus T $ and an ordered $k$-tuple of planar reduced
trees that represent operations in $\Mag^S$. Thus we get the
formula $f_{\magroot{T}{S}}(x)=x+\sum_{k\in S\setminus
T}(f_{\Mag^S}(x))^k.$
\end{proof}

\begin{proposition}\label{prop:combinatorics}
Let $T\subset S\subset \NN$. Then:

$$x+\sum_{k\in S\setminus
T}(f_{\Mag^S}(x))^k=f_{\magroot{T}{S}}(x)=(x - \sum_{i\in T}
x^i)\circ (x - \sum_{i\in S} x^i)^{\circ(-1)}.$$
\end{proposition}

\begin{proof}
The first equality holds by Lemma \ref{lem:combinatorics}. For the
second equality, we may use  Corollary \ref{corollary:pbw}
together with Proposition \ref{prop:inversion}:

It follows from  Corollary \ref{corollary:pbw} that
$$(f_{\Mag^T}(x))^{\circ(-1)}\circ
f_{\Mag^S}(x)=f_{\magroot{T}{S}}(x).$$

By Proposition \ref{prop:inversion}, $(f_{\Mag^T}(x))^{\circ(-1)}$
is given by $(x - \sum_{i\in T} x^i)$, and also $f_{\Mag^S}(x)$ is
the composition inverse of $(x - \sum_{i\in S} x^i)$. Hence the
second equation follows.

\end{proof}

\begin{example}

The Super-Catalan numbers (cf.\ \cite{S} A001003) are the
coefficients of the generating series of $\Mag^{\NN}$:
\begin{equation*}
f^{\Mag^{\NN}}(x) =\sum_{n\geq 1}C_n x^n=
\frac{1}{4}(1+x-\sqrt{1-6x+x^2}).
\end{equation*}

Let $S=\NN$ and $T=\NN-\{2\}=\{3,4,\ldots\}$. Then we obtain

$$f_{\magroot{T}{S}}(x)=x+\sum_{k\in \{2\}}(f_{\Mag^{\NN}}(x))^k=
x+\frac{(x-1)^2-(x+1)\sqrt{1-6x+x^2}}{8}.$$

In low degrees, this series is equal to $x+x^2+2x^3+7x^4+ 28x^5+
121x^6+ 550 x^7+\ldots$

To verify the second equation in this case, we note that
\begin{eqnarray*}
&(x - \sum_{i\in \NN} x^i)^{\circ(-1)}\\
=&(\frac{x(1-2x)}{1-x})^{\circ(-1)}\\
=&\frac{1}{4}(1+x-\sqrt{1-6x+x^2})\\ =& f_{\Mag^{\NN}}(x).
\end{eqnarray*}

Since $(x - \sum_{i\geq 3} x^i)\circ f_{\Mag^{\NN}}(x) =
f_{\Mag^{\NN}}(x) - \sum_{i\geq 3}(f_{\Mag^{\NN}}(x))^i,$ we get
the equation

$$f_{\Mag^{\NN}}(x)=f_{\magroot{T}{S}}(x)+ \sum_{i\geq
3}(f_{\Mag^{\NN}}(x))^i =x+\sum_{k\geq 2}(f_{\Mag^{\NN}}(x))^k$$

in this case.
\end{example}

\begin{example}

Let $S=\NN$ and $T=\{2k + 1 : k\geq 1\}$. Then we obtain by
analogous computations that

$$f_{\magroot{T}{S}}(x)=x-\frac{1+x^2-2x
-(x+1)\sqrt{1-6x+x^2}}{(-7)+x^2-2x -(x+1)\sqrt{1-6x+x^2}}.$$

In low degrees, this series is equal to $x+x^2+2x^3+8x^4+ 32x^5+
140x^6+ 640 x^7+\ldots$

\end{example}

\bigskip

 \noindent {\bf Acknowledgement.}

 R.~Holtkamp wants to thank for the hospitality of the Institut de Math\'ematiques et de mod\'elisation de
 Montpellier.
Both authors thank the Max-Planck-Institute for Mathematics Bonn,
where this work was begun, for support. The authors thank
Jean-Louis Loday for helpful discussions and his remarks on a
first version of this paper.
Both authors would like to thank the referee for helpful comments.


\begin{thebibliography}{999}
\bibitem{BM} D.\ Borisov, Y.\ Manin, Generalized operads and their inner
cohomomorphisms, Preprint math.CT/0609748.
\bibitem{B1}  E.\ Burgunder, Big\`ebre magmatique et big\`ebre Associative-Zinbiel, M\'emoire de Master, Strasbourg 2005.
\bibitem{B2}  E.\ Burgunder, Infinite magmatic bialgebras,
                 \textsl{Adv.\ Appl.\ Math.} 40 (2008), no. 3, 309-329.
\bibitem{Ca} P.\ Cartier, Hyperalg\`ebres et groupes de Lie formels, S\'eminaire Sophus Lie, 2e ann\'ee: 1955/56.
 Facult\'e des Sciences de Paris.
\bibitem{G} L.\ Gerritzen, Planar rooted trees and non-associative exponential series,
 \textsl{Adv.\ Appl.\ Math.} 33 (2004), no. 2, 342--365.
\bibitem{GK} V.\ Ginzburg, and M.\ Kapranov, Koszul duality for operads, \textsl{Duke Math.\ J.}
76 (1994), no. 1, 203--272.
\bibitem{H} R.\ Holtkamp, On Hopf algebra structures over free operads, \textsl{Adv.\ Math.} 207 (2006), 544--565.
\bibitem{HLR} R.\ Holtkamp, J.-L.\ Loday, M.\ Ronco, Coassociative magmatic bialgebras and the Fine numbers,
 Preprint math.RA/0609125, to appear in {J.\ Algebraic Combin.}
\bibitem{L} J.-L.\ Loday, Generalized bialgebras and triples of operads, Preprint math.QA/0611885.
\bibitem{LR1}  J.-L.\ Loday, M.\ Ronco,  On the structure of cofree Hopf algebras,
\textsl{J.\ Reine Angew.\ Math.} 592 (2006), 123--155.
\bibitem{LR2} J.-L.\ Loday,  M.\ Ronco, Alg\`ebres de Hopf colibres,
\textsl{C.\ R.\ Acad.\ Sci.\ Paris} 337 (2003), no. 3, 153--158.
\bibitem{M}  J.\ Milnor, J.\ C.\ Moore, On the structure of Hopf
algebras, \textsl{Ann.\ of Math.} (2) 81 (1965), 211--264.
\bibitem{Q}  D.\ Quillen, Rational homotopy theory, \textsl{Ann.\ of Math.}(2) 90 (1969), 205--295.
\bibitem{R1} M.\ Ronco,  Eulerian idempotents and Milnor-Moore theorem for certain non-commutative Hopf algebras,
\textsl{J. Algebra} 254 (2002), no. 1, 152--172.
\bibitem{R}  M.\ Ronco, Primitive elements in a free dendriform algebra. New trends in Hopf algebra theory (La Falda,
1999), 245--263, \textsl{Contemp.\ Math.} 267, Amer.\ Math.\ Soc.,
Providence, RI, 2000.
\bibitem{S} N.\ Sloane (Edt.), The On-Line Encyclopedia of Integer
                      Sequences, 2007,
                       \newline
                       http://www.research.att.com/\~\ \!\!njas/sequences/index.html
\bibitem{V}  B.\ Vallette,  A Koszul duality for props, \textsl{Trans.\ Amer.\ Math.\ Soc.} 359 (2007), 4865--4943.
\end{thebibliography}
\end{document}